\documentclass[a4paper,10pt]{article}
\usepackage{amsmath}
\usepackage{amssymb}
\usepackage{amsthm}
\usepackage{mathrsfs}
\usepackage{mathtools}
\usepackage{bm}
\usepackage{fullpage}
\usepackage{enumerate}
\usepackage{hyperref}
\usepackage{slashed}
\usepackage{marvosym}
\usepackage{caption}
\usepackage{graphicx}

\newcommand{\upd}{\mathrm{d}}
\newcommand{\D}{\mathscr{D}}
\newcommand{\Lbar}{\underline{L}}

\newcommand{\dVol}{\mathrm{d}\textit{vol}}

\allowdisplaybreaks[1]

\usepackage[utf8]{inputenc}
\usepackage[backend=biber, style=numeric-comp, sorting=none]{biblatex}
\renewbibmacro*{in:}{%
	\ifentrytype{article}{}{\printtext{\bibstring{in}\intitlepunct}}}
\newtheorem{theorem}{Theorem}[section]

\newtheorem{corollary}[theorem]{Corollary}

\theoremstyle{definition}

\newtheorem{condition}[theorem]{Condition}

\theoremstyle{remark}

\begin{document}

\title{Global existence for systems of nonlinear wave equations with bounded, stable asymptotic systems}
\author{Joe Keir \\ \\
{\small DAMTP, Centre for Mathematical Sciences, University of Cambridge}, \\
{\small \sl Wilberforce Road, Cambridge CB3 0WA, UK} \\ \\
{\small Trinity College, Cambridge} \\ \\
\small{j.keir@damtp.cam.ac.uk}}
\maketitle

\begin{abstract}
	Some systems of nonlinear wave equations admit global solutions for all sufficiently small initial data, while others do not. The (classical) null condition guarantees that such a result holds, but it is too strong to capture certain systems -- most famously the Einstein equations -- which nevertheless admit global solutions for small initial data. The \emph{weak null condition} has been proposed as a sufficient condition for such a result to hold; it takes the form of a condition on a related set of nonlinear ODEs known as the ``asymptotic system''. Previous results in this direction have required certain structural conditions on the asymptotic system in addition to the weak null condition. In this work we show that, if the solutions to the asymptotic system are bounded (given small initial data), and, in addition, if these solutions are stable against rapidly decaying perturbations, then the corresponding system of nonlinear wave equations admits global solutions for all sufficiently small initial data. This avoids any direct assumptions on the structure of the nonlinear terms. We also give an example of a class of systems obeying this condition but not obeying previously identified structural conditions. For this class the asymptotic system arises as a generalisation of the ``Euler equations'' for rigid body motion, associated with a left-invariant Hamiltonian flow on a finite dimensional Lie group. This work relies heavily on the author's previous work \cite{Keir2018}.
\end{abstract}

\section{Introduction}

In three spatial dimensions, the global behaviour of systems nonlinear wave equations can differ drastically from one system to another, even in the case of small initial data. For example, the scalar wave equation
\begin{equation}
\label{equation null condition example}
	\Box \phi = (\partial_t\phi)^2 - \sum_{i = 1}^3 (\partial_i \phi)^2
\end{equation}
admits global solutions for all sufficiently small initial data \cite{Klainerman1986, Christodoulou1986}, whereas \emph{all} non-trivial solutions to the superficially similar wave equation
\begin{equation}
\label{equation John's example}
	\Box \phi = (\partial_t\phi)^2
\end{equation}
blow up in finite time \cite{John1981}.

These systems can be distinguished using the ``null condition'' \cite{Klainerman1980}: the first system satisfies the null condition, whereas the second does not. Indeed, Klainerman and (independently) Christodoulou proved that all systems obeying the null condition admit global solutions for sufficiently small initial data \cite{Klainerman1986, Christodoulou1986}. The null condition therefore serves as a sufficient condition for a system of nonlinear wave equations to admit global solutions for small initial data.

However, the null condition not a \emph{necessary} condition for a small data global existence result to hold. This is illustrated by the system
\begin{equation}
\label{equation weak null condition example}
\begin{split}
	\Box \phi_1 &= 0 \\
	\Box \phi_2 &= (\partial_t \phi_1)^2
\end{split}
\end{equation}
which does \emph{not} satisfy the null condition, but which nevertheless admits global solutions for all small initial data. A more interesting, and much more complicated example is provided by the vacuum Einstein equations in wave coordinates, which also fails to satisfy the null condition but which nevertheless admits global solutions given small initial data \cite{Lindblad2004, Lindblad2005}.

As an attempt to classify such systems, we can examine the associated \emph{asymptotic systems} \cite{Hoermander1987, Hoermander1997}. This is a system of nonlinear, first order ODEs, obtained from the corresponding system of nonlinear wave equations by a process that involves dropping all of the terms that are expected to exhibit ``good'' behaviour\footnote{See, for example, \cite{Keir2018} and references therein for further details. Note that in \cite{Keir2018} we provide a prescription for handling quasilinear systems, but for simplicity here we will only discuss semilinear equations.}. For example, the asymptotic system associated to the wave equation \eqref{equation null condition example} is simply\footnote{In these equations, the parameter $s$ can be identified with $\log (1+r)$ along outgoing null rays.} the linear ODE
\begin{equation*}
	\frac{\partial}{\partial s} \Phi = 0
\end{equation*}
which clearly has global (bounded) solutions, whereas the asymptotic system for the equation \ref{equation John's example} is
\begin{equation*}
	\frac{\partial}{\partial s} \Phi = \frac{1}{4} \Phi^2
\end{equation*}
which admits solutions with arbitrarily small initial data, but which blow up at a finite value of the parameter $s$. Finally, the asymptotic system for the equation \eqref{equation weak null condition example} is
\begin{equation*}
\begin{split}
	\frac{\partial}{\partial s} \Phi_0 &= 0
	\\
	\frac{\partial}{\partial s} \Phi_1 &= \frac{1}{4} (\Phi_0)^2
\end{split}
\end{equation*}
Although this is a nonlinear system, it is easy to see that it admits global solutions, and moreover these solutions grow like $\epsilon^2 s$ for initial data of size $\epsilon$.

In general, a system of nonlinear wave equations is said to obey the \emph{weak null condition} if the corresponding asymptotic system admits global solutions, \emph{and} if those solutions grow no faster than $e^{\epsilon s}$, where $\epsilon$ measures the size of the initial data \cite{Lindblad2003}. For example, the asymptotic system associated with the system of wave equations
\begin{equation*}
\begin{split}
	\Box \phi_1 &= 0 \\
	\Box \phi_2 &= (\partial_t \phi_1)(\partial_t \phi_2) \\
	\Box \phi_3 &= (\partial_t \phi_2)(\partial_t \phi_3)
\end{split}
\end{equation*}
does admit global solutions, but these solutions grow at the super-exponential rate $e^{e^{\epsilon s}}$. Hence this system does not obey the weak null condition.

It has been suggested that, similarly to the classical null condition, the weak null condition provides a sufficient condition for a system of wave equations to admit global solutions for small initial data. Indeed, several examples of wave equations satisfying the weak null condition have been studied in detail, and in each case a small-data global-existence result was proved for the example systems. Proving such a result for the system \eqref{equation weak null condition example} is trivial\footnote{In fact, we do not even require small initial data for this system!}; more complicated semilinear systems were studied in \cite{Alinhac2006}, a quasilinear scalar wave equation was studied in \cite{Lindblad2008} and the Einstein equations in wave coordinates were studied in \cite{Lindblad2004, Lindblad2005}.

Recently, we proved a small-data global-existence result for a large class of nonlinear wave equations, which includes all of the examples mentioned above as special cases. This result was proven under a specific structural condition on the semilinear terms, which we referred to as the \emph{hierarchical weak null condition}. Glimpses of this structure can be seen in the example systems above: essentially, we require that the ``bad'' semilinear terms appearing in the equation for a field at some given level of the hierarchy only involve certain combinations of fields at lower levels of the hierarchy.

This structural condition guarantees that the weak null condition is satisfied, however, as we shall see in this paper, there are examples of wave equations satisfying the weak null condition but which fail to satisfy the hierarchical null condition. It is important to note that all of the previously studied examples of systems obeying the weak null condition (that we are aware of) in fact do obey the hierarchical condition.

In this paper we will consider a different class of nonlinear wave equations. The associated asymptotic systems will not, in general, exhibit the hierarchical structure required to satisfy the hierarchical null condition. Instead, the asymptotic system will be required to obey the following two conditions:
\begin{enumerate}
	\item All solutions of the asymptotic system arising from small initial data are \emph{bounded}.
	\item Solutions to the asymptotic system are stable to small, exponentially decaying perturbations (see condition \ref{Condition boundedness and stability} for precise details).
\end{enumerate}
Note that, as in our previous work \cite{Keir2018}, this condition refers \emph{only} to the asymptotic system, and not to the original system of nonlinear wave equations. However, unlike this previous work, the conditions above are stated only in terms of the \emph{behaviour of solutions} to the asymptotic system -- we do not require any explicit structure in the asymptotic system itself. As such, these conditions are closer in spirit to the original weak null condition.

We will also give an example of a class of asymptotic systems obeying the conditions above. These systems obey both the boundedness and stability criteria due to the presence of a Hamiltonian structure in the asymptotic system, and might be of independent interest. In fact, the asymptotic system will correspond to the flow generated on the cotangent space at the origin associated with a left-invariant Hamiltonian on some finite dimensional Lie group, where the Hamiltonian vanishes quadratically at the origin. These systems are generalisations of the Euler equations for the motion of rigid bodies -- see \cite{Mienko1978}.

Concrete example systems have played a prominent role in the development of the theory of global existence for nonlinear waves, as illustrated by our discussion above. With this in mind, we give here a model system that falls into the class of example systems discussed above, i.e.\ its asymptotic system arises from a Hamiltonian on a Lie group. This system consists of a set of three coupled semilinear scalar waves, obeying the equations
\begin{equation}
\begin{split}
\label{equation model Euler equations}
	\Box \phi_{(1)} &= \frac{ I_{(2)} - I_{(3)} }{ 4I_{(1)} } (\partial_t \phi_{(2)}) (\partial_t \phi_{(3)})
	\\
	\Box \phi_{(2)} &= \frac{ I_{(3)} - I_{(1)} }{ 4I_{(2)} } (\partial_t \phi_{(3)}) (\partial_t \phi_{(1)})
	\\
	\Box \phi_{(3)} &= \frac{ I_{(1)} - I_{(2)} }{ 4I_{(3)} } (\partial_t \phi_{(1)}) (\partial_t \phi_{(2)})
\end{split}
\end{equation}
where $I_{(1)}$, $I_{(2)}$ and $I_{(3)}$ are positive constants. The corresponding asymptotic system is the system of ODEs
\begin{equation}
\begin{split}
	I_{(1)} \dot{\Phi}_{(1)} + \left( I_{(2)} - I_{(3)} \right) \Phi_{(2)} \Phi_{(3)}
	\\
	I_{(2)} \dot{\Phi}_{(2)} + \left( I_{(3)} - I_{(1)} \right) \Phi_{(3)} \Phi_{(1)}
	\\
	I_{(3)} \dot{\Phi}_{(3)} + \left( I_{(1)} - I_{(2)} \right) \Phi_{(1)} \Phi_{(2)}
\end{split}
\end{equation}
These are the Euler equations for rigid body motion, with the constants $I_{(a)}$ corresponding to the principle moments of inertia and the $\Phi_{(a)}$ corresponding to the angular velocity of a rigid body around the corresponding principle axis. Note that the solutions $\Phi_{(a)}$ are bounded; in fact, if we define
\begin{equation*}
	H(\bm{\Phi}) := \frac{1}{2} \sum_{a = 1}^3 I_{(a)} \Phi_{(a)}^2
\end{equation*}
then $H$ is constant on the flow generated by the asymptotic system. For the corresponding wave equations we find that, along any outgoing null ray, $H( r(\partial_t - \partial_r)\bm{\phi} )$ approaches a constant value as $r \rightarrow \infty$.

\begin{figure}[phtb]
	\centering
	\includegraphics[height = 0.6\textheight, keepaspectratio]{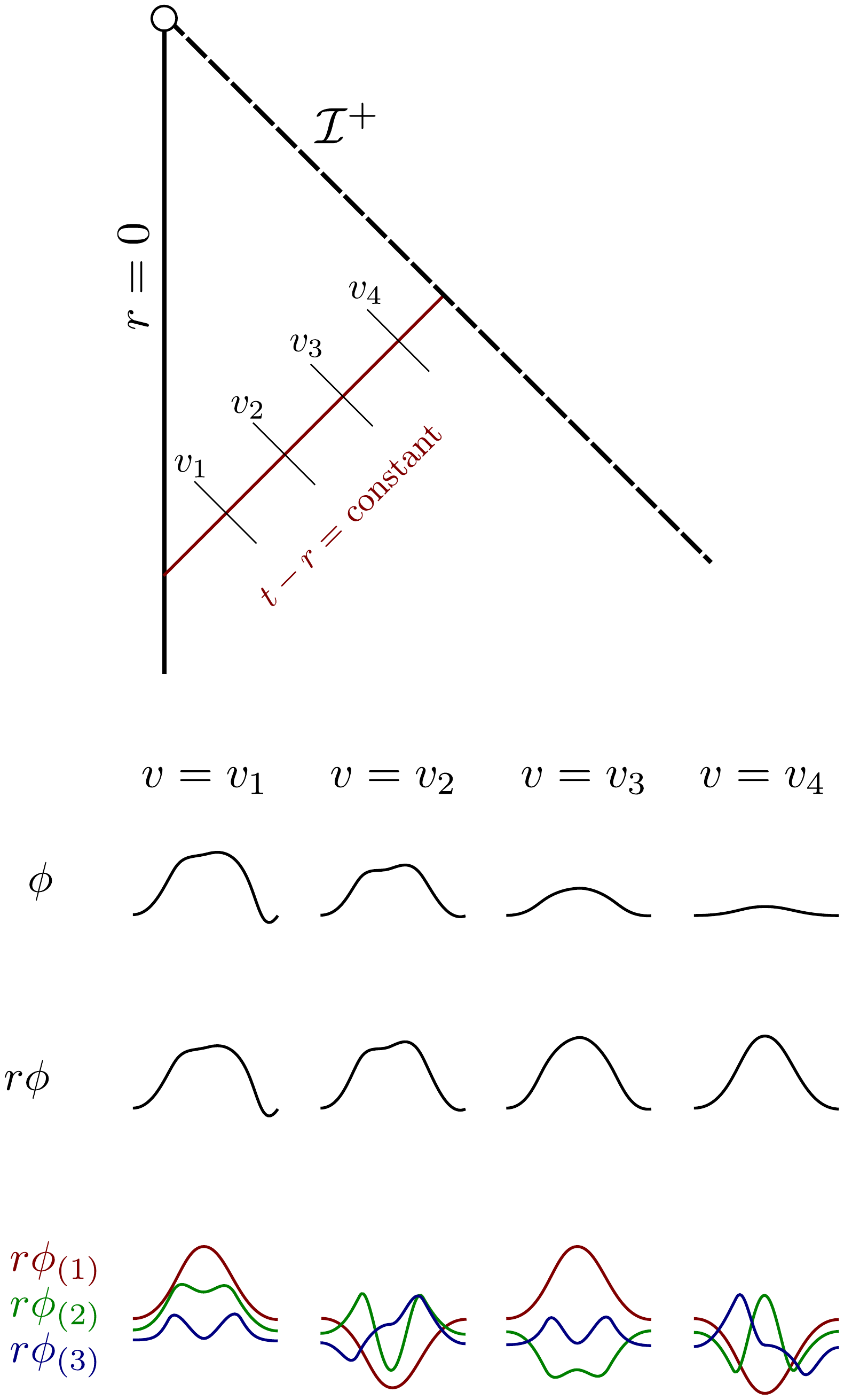}
	\caption[figure1]{
		A figure illustrating the behaviour of various fields towards null infinity. At the top we have sketched a Penrose diagram of a portion of Minkowski space. We are interested in the behaviour of various fields near to the outgoing null cone $u = (t-r) = \text{constant}$, shown in red on the diagram. We have also sketched various segments of the ingoing null cones $v = (t+r) = v_1$, $v_2$, $v_3$ and $v_4$, near to their intersections with this outgoing null cone.
		\\
		Below the Penrose diagram we have illustrated the behaviour of various fields on the segments of the ingoing null cones mentioned above. In reality, the values of these fields will also depend on the angular coordinates, but for illustrative purposes we have shown the fields at some fixed values of the angular coordinates. These fields can then be thought of as functions of $u = (t-r)$ at the corresponding fixed values of $v = (t+r)$, and so the $x$-axis corresponds to $u$.
		\\
		In the first row we have shown the value of the field $\phi$, which solves either a linear wave equation or a nonlinear wave equation with the classical null condition. As $v \rightarrow \infty$, we have $\phi \rightarrow 0$. This is also true for fields satisfying wave equations which only obey various versions of the weak null condition.
		\\
		In the second row we have shown the value of the ``radiation field'' $r\phi$. For each angle, and for each value of $u$, this tends to some finite limit as $v \rightarrow \infty$. In other words, there is a function $\tilde{\psi}$ of $u$ and the angular coordinates such that $(r\phi - \tilde{\psi}) \rightarrow 0$.
		\\
		In the third row we have shown the three ``radiation fields'' $r\phi_{(i)}$ from the model system \eqref{equation model Euler equations}. These fields are bounded as $v \rightarrow \infty$, but along some given null line (e.g.\ at fixed $u$ and for fixed angular coordinates) the three fields continue to evolve as $v \rightarrow \infty$, albeit at a (logarithmically) decreasing rate. For large values of $v$, the radiation fields asymptotically obey the Euler equations, with $\log v$ (or $\log r$) playing the role of the time parameter.
	}
	\label{figure Euler Equations F1}
\end{figure}

We note here that the scattering theory for such equations promises to be extremely interesting. It is likely the case that, along some fixed null ray, the ``radiation field'' $r \bm{\phi}$ (where $\bm{\phi}$ is the vector of solutions to the system of nonlinear wave equations) behaves similarly to the solutions to the corresponding asymptotic system. For linear wave equations, or for equations satisfying the classical null condition, this means that $r\bm{\phi}$ will tend to some constant vector along each fixed null ray, which is indeed the case. On the other hand, for the types of equation considered in this paper, $r\bm{\phi}$ is bounded but might not have a limit along a given null ray. In fact, for solutions arising from Hamiltonians in the manner described above, along each null ray we expect $r\bm{\phi}$ to asymptotically approach a solution to the corresponding dynamical system. See figure \ref{figure Euler Equations F1} for a detailed picture of this scenario.

We also note that it would be possible to generalise our construction further, to take the kinds of systems studied in this paper and to place a hierarchical system (i.e.\ one satisfying the ``hierarchical null condition'' of \cite{Keir2018}) ``on top'' of this system. There may be other ways to generalise our construction, but we do not pursue them here.

The rest of this paper is organised as follows. In section \ref{section notation and theorem} we will state the precise form of the theorem that we prove, and define the notation used in the rest of the paper. Then, instead of reproving the global existence result of \cite{Keir2018} from scratch, in section \ref{section proving theorem} we will instead point out the places in \cite{Keir2018} where the hierarchical null condition is used, and show that, in each case, this condition can be replaced with the boundedness and stability conditions. Then, in section \ref{section Euler equations} we will show that these conditions are satisfied by all systems of wave equations where the asymptotic system arises from a Hamiltonian in the way described above.

\section*{Acknowledgements}
I am very grateful to Igor Rodnianski for his interest in this problem, and for suggesting systems arising from Hamiltonians as a model. I am also grateful for the hospitality of Princeton University, where some of this work was done.

\section{Notation and a precise statement of the theorem}
\label{section notation and theorem}
 
Consider a system of nonlinear wave equations
\begin{equation*}
	\Box_g \phi = F(\phi, \partial \phi)
\end{equation*}
where $g = g(\phi)$ is a Lorentzian metric, and where
\begin{equation*}
	\phi : \mathbb{R}^4 \rightarrow \mathcal{V}
\end{equation*}
for some finite-dimensional vector space $\mathcal{V}$. Here also, $F$ is a nonlinear term of the form
\begin{equation*}
	F(\phi, \partial \phi) = F^{(0)} + F^{(1)}(\partial \phi, \partial \phi) + F^{(2)}(\phi, \partial \phi)
\end{equation*}
where $F^{(0)}$ is some homogeneous term satisfying suitable bounds, $F^{(1)}$ is a quadratic form and $F^{(2)}$ vanishes at least cubically at $(\phi, \partial\phi) = 0$. Suppose also that the metric $g$ can be written, in some coordinate system $(x^a)$, as
\begin{equation*}
	g_{ab} = m_{ab} + h_{ab}(\phi)
\end{equation*}
where $m_{ab} = \text{diag}(-1,1,1,1)$ is the Minkowski metric in inertial coordinates. Define $t = x^0$ and $r^2 = \sum_{i = 1}^3 (x^i)^2$.

Define $u$ as an outgoing solution of the eikonal equation
\begin{equation*}
\begin{split}
	g^{-1}(\upd u, \upd u) &= 0 \\
	u\big|_{r = r_0} &= t - r_0
\end{split}
\end{equation*}
and define the (future directed, null) vector field $L$ as
\begin{equation*}
	L_{(\text{geo})} = - (\upd u)^\sharp
\end{equation*}
Now define the normalised null vector field
\begin{equation*}
	L = \frac{1}{L_{(\text{geo})}(r)} L_{(\text{geo})}
\end{equation*}
so that $L(r) = 1$. Define the conjugate null vector field $\Lbar$ as follows:
\begin{itemize}
	\item $\Lbar$ is null: $g(\Lbar,\Lbar) = 0$
	\item $\Lbar$ is orthogonal (with respect to $g$) to the ``spheres'' of constant $u$ and $r$
	\item $\Lbar$ is normalised by $g(L, \Lbar) = -2$
\end{itemize}
Next, define the connection coefficient $\omega$ by
\begin{equation*}
	\omega = -\frac{1}{2} g(\D_L L, \Lbar)
\end{equation*}
where $\D$ is the covariant derivative with respect to $g$. Finally, we define $\slashed{\nabla}$ as the restriction of the covariant derivative $\D$ to the tangent space of the spheres of constant $u$ and $r$. That is, $\slashed{\nabla}_X \phi = \D_X \phi$ if the vector field $X$ is tangent to the spheres, while $\slashed{\nabla}_L \phi = \slashed{\nabla}_{\Lbar} \phi = 0$.

Now, we can define the ``reduced wave operator'': for any scalar field $\phi$, we have
\begin{equation*}
	\tilde{\Box}_g \phi := \Box_g \phi + \omega (\Lbar \phi)
\end{equation*}
Note that $\omega$ can be written, essentially, in terms of derivatives of $\phi$ (see \cite{Keir2018} for details).

Using the reduced wave operator, we can write our system of nonlinear wave equations as follows:
\begin{equation}
\label{equation system of nonlinear reduced wave equations}
	\tilde{\Box} \phi_{(A)}
	=
	F_{(A)}^{(BC)} (\Lbar \phi_{(B)}) (\Lbar \phi_{(C)}) + \text{good terms}
\end{equation}
where the $F_{(A)}^{(BC)}$ are constants, and where the ``good terms'' are either homogeneous terms that are independent of $\phi$ and which satisfy suitable bounds, higher order terms (i.e.\ cubic in the fields $\phi_{(A)}$ and their derivatives), or quadratic terms involving a pair of derivatives of $\phi$, where at least one derivative is in a direction tangent to the level sets of $u$.

The asymptotic system corresponding to this set of wave equations \eqref{equation system of nonlinear reduced wave equations} is the following system of nonlinear ODEs:
\begin{equation}
\label{equation asymptotic system}
	\frac{\partial}{\partial s} \Phi_{(A)}
	=
	-F_{(A)}^{(BC)} \Phi_{(B)} \Phi_{(C)}
\end{equation}

Now we can make precise the conditions we require of the asymptotic system:
\begin{condition}[Boundedness and stability for the asymptotic system]
\label{Condition boundedness and stability}

	Given some small $\delta > 0$ and some constant $C > 0$, there exists some $\epsilon > 0$ and some constant $\tilde{C} > 0$ (depending on $C$ and $\delta$) such that, if $f_{(A)}(s)$ are any set of functions satisfying the bound
		\begin{equation*}
			|f_{(A)}(s)| < C \epsilon e^{-\delta s}
		\end{equation*}
		then all solutions to the system
		\begin{equation}
			\frac{\partial}{\partial s} \Phi_{(A)}
			=
			F_{(A)}^{(BC)} \Phi_{(B)} \Phi_{(C)}
			+ f_{(A)}(s)
		\end{equation}
		arising from initial data satisfying
		\begin{equation*}
			\sup_{(A)} \left( \Phi_{(A)} \big|_{s = 0} \right) < \epsilon
		\end{equation*}
		obey the bound
		\begin{equation*}
			\sup_{s \geq 0} \sup_{(A)} \left( \Phi_{(A)}(s)  \right) < \tilde{C}\epsilon
		\end{equation*}
\end{condition}

\vspace{3mm}
Note that, if a system satisfies condition \ref{Condition boundedness and stability}, then, by choosing $f(s) = 0$, we find that the solutions to the asymptotic system \eqref{equation asymptotic system} with small initial data are uniformly bounded by some small constant.

We are now in a position to state the main theorem:

\begin{theorem}[Global existence for small initial data for systems of wave equations with asymptotic systems obeying condition \ref{Condition boundedness and stability}]
	\label{theorem main theorem}
	
	Let $\phi : \mathcal{M} \rightarrow \mathcal{V}$, where $\mathcal{V}$ is some finite-dimensional vector space, and suppose that $\phi$ satisfies the equation
	\begin{equation*}
		\tilde{\Box}_g \phi = F^{(0)} + F^{(1)}(\partial \phi, \partial \phi) + F^{(2)}(\phi, \partial \phi)
	\end{equation*}
	
	Suppose that the asymptotic system associated with this set of nonlinear wave equations obeys the condition \ref{Condition boundedness and stability}, i.e.\ it has stably bounded solutions.
	
	We require the following bounds for the inhomogeneous terms which are independent of the field $\bm{\phi}$ and its derivatives. Define $F^{(0)}_{(A, n)} = \mathscr{Y}^n F_{(A)}^{(0)}$, where each operator $\mathscr{Y}$ is any operator from the set $\{ \slashed{\D}_T , r\slashed{\nabla}, r\slashed{\D}_L \}$ (see \cite{Keir2018} for this notation). Then we require the following bounds for these inhomogeneous terms $F_{(A)}^{(0)}$:
	
	Pointwise bounds: for all $n \leq N_1$, we have
	\begin{equation*}
	|F^{(0)}_{(A,n)}| \leq \epsilon^{(N_2 - 2 - n)} (1+r)^{-2 + 2C_{(n)}\epsilon} (1+\tau)^{-\beta}
	\end{equation*}
	for some sufficiently large constants $C_{(n)}$, which satisfy
	\begin{equation*}
	C_{(n)} \gg C_{(n-1)}
	\end{equation*}
	and
	\begin{equation*}
	C_{(0)} = 0
	\end{equation*}
	
	We also require some $L^2$ bounds on these homogeneous terms. In order to state these, we first decompose the homogeneous terms as follows:
	\begin{equation*}
	\begin{split}
	F_{(A,n)}^{(0)}
	&=
	F_{(A,n,1)}^{(0)}
	+ F_{(A,n,2)}^{(0)}
	+ F_{(A,n,3)}^{(0)}
	\\
	&=
	F_{(A,n,4)}^{(0)}
	+ F_{(A,n,5)}^{(0)}
	+ F_{(A,n,6)}^{(0)}
	\end{split}	
	\end{equation*}
	then we require the $L^2$ bounds	
	\begin{equation*}
	\begin{split}
	&\int_{\mathcal{M}_{\tau_0}^\tau} \epsilon^{-1} \bigg(
	(1+r)^{1-C_{[n]}\epsilon}|F^{(0)}_{(A,n)}|^2
	+ (1+r)^{\frac{1}{2}\delta}(1+\tau)^{1+\delta} |F^{(0)}_{(A,n,1)}|^2
	+ (1+r)^{1-3\delta}(1+\tau)^{2\beta} |F^{(0)}_{(A,n,2)}|^2
	\\
	&\phantom{\int_{\mathcal{M}_{\tau_0}^\tau} \epsilon^{-1} \bigg(}
	+ (1+r)^{1+\frac{1}{2}\delta} |F^{(0)}_{(A,n,3)}|^2
	\bigg)\dVol_g
	\lesssim
	\frac{1}{C_{[n]}} \epsilon^{2(N_2 - n + 2)}(1+\tau)^{-1+C_{(n)}\delta}
	\\
	\\
	&\int_{\mathcal{M}_{\tau_0}^\tau \cap \{r \geq r_0\}} \epsilon^{-1}\bigg(
	r^{1-C_{[n]}\epsilon} (1+\tau)^{1+\delta} |F^{(0)}_{(A,n,4)}|^2
	+ r^{2-C_{[n]}\epsilon-2\delta} (1+\tau)^{2\beta} |F^{(0)}_{(A,n,5)}|^2
	\\
	&\phantom{\int_{\mathcal{M}_{\tau_0}^\tau \cap \{r \geq r_0\}} \epsilon^{-1}\bigg(}
	+ r^{2-C_{[n]}\epsilon} |F^{(0)}_{(A,n,6)}|^2
	\bigg)\dVol_g
	\lesssim
	\frac{1}{C_{[n]}} \epsilon^{2(N_2 - n + 2)}(1+\tau)^{C_{(n)}\delta}
	\end{split}
	\end{equation*}
	where $\delta$ is some sufficiently small constant and the $C_{[n]}$ are some sufficiently large constants, satisfying
	\begin{equation*}
	C_{[n]} \gg C_{[n-1]}
	\end{equation*}
	and also, for all $n_1$, $n_2$
	\begin{equation*}
	C_{[n_1]} \gg C_{(n_2)}
	\end{equation*}

	Next, suppose that the rectangular components of the metric components can be expressed as
	\begin{equation*}
	h_{ab}(x, \bm{\phi}) = h^{(0)}_{ab}(x) + h^{(1)}_{ab}(x, \bm{\phi})
	\end{equation*}
	such that the following bounds hold: for all $n \leq N_1$
	\begin{equation*}
	\begin{split}
	|\mathscr{Y}^n h^{(1)}_{ab}| 
	&\lesssim
	\sum_{m \leq n} \sum_{(a)} |\mathscr{Y}^m \phi_{(a)}| + \mathcal{O}(|\mathscr{Y}^{\leq n} \bm{\phi}|^2)
	\\
	|\slashed{\D} \mathscr{Y}^n h^{(1)}_{ab}| 
	&\lesssim
	\sum_{m \leq n} \sum_{(a)} |\slashed{\D} \mathscr{Y}^m \phi_{(a)}| + \mathcal{O}\left( \sum_{j+k \leq n}|\slashed{\D}\mathscr{Y}^{j} \bm{\phi}||\mathscr{Y}^{k} \bm{\phi}| \right)
	\\
	|\overline{\slashed{\D}} \mathscr{Y}^n h^{(1)}_{ab}| 
	&\lesssim
	\sum_{m \leq n} \sum_{(a)} |\overline{\slashed{\D}} \mathscr{Y}^m \phi_{(a)}| + \mathcal{O}\left( \sum_{j+k \leq n}|\slashed{\D}\mathscr{Y}^{j} \bm{\phi}||\mathscr{Y}^{k} \bm{\phi}| \right)
	\end{split}
	\end{equation*}
	(see \cite{Keir2018} for the notation used here) and also such that we have the following bound
	\begin{equation*}
	|\partial h^{(1)}|_{LL} = |\partial h^{(1)}_{ab}|L^a L^b \lesssim \sum_{(A) \, | \, \phi_{(A)} \in \Phi_{[0]}} |\partial \phi|_{(A)} + \mathcal{O}(|\bm{\phi}||\partial \bm{\phi}|)
	\end{equation*}
	Additionally, the lower order terms in the metric perturbations are required to satisfy the following pointwise bounds, for all $n \leq N_1$
	\begin{equation*}
	\begin{split}
	|\mathscr{Y}^n h^{(0)}_{ab}| &\leq \frac{1}{2}\epsilon (1+r)^{-\frac{1}{2} + \delta}
	\\
	|\slashed{\D}\mathscr{Y}^n h^{(0)}_{ab}| &\leq \begin{cases}
	\frac{1}{2}\epsilon \left( (1+r)^{-1} + (1+r)^{-1+\delta} (1+\tau)^{-\beta} \right) \\
	\frac{1}{2}\epsilon (1+r)^{-1+C_{(n)}\epsilon}
	\end{cases}
	\\
	|\slashed{\D}\mathscr{Y}^n h^{(0)}_{ab}| &\leq \frac{1}{2}\epsilon (1+r)^{-1-\delta}
	\end{split}
	\end{equation*}

	We also require the following bounds, giving additional control over lower order terms:
	\begin{equation*}
	\begin{split}
	|\partial h^{(0)}_{ab}| &\leq \frac{1}{2}\epsilon (1+r)^{-1 + C\epsilon}(1+\tau)^{-C^* \delta}
	\\
	|\slashed{\D} \mathscr{Y} h^{(0)}_{ab}| &\leq \frac{1}{2}\epsilon (1+r)^{-1 + C_{(1)}\epsilon}(1+\tau)^{-C^* \delta} 
	\\
	|\partial h^{(0)}_{ab}|L^a L^b &\leq \frac{1}{2} \epsilon (1+r)^{-1}
	\end{split}
	\end{equation*}

	Finally, we suppose that the initial data for the fields $\phi_{(a)}$ is posed on the hypersurface $\Sigma_{\tau_0}$, which consists of two parts: the hypersurface\footnote{Here, as elsewhere, $r$ is defined relative to the \emph{rectangular coordinates} by $r = \sqrt{(x^1)^2 + (x^2)^2 + (x^3)^2}$.} $\{ x^0 = t_0 = \text{constant}\} \cap \{r \leq r_0\}$ together with an outgoing \emph{characteristic hypersurface} emanating from the sphere $r = r_1$, $t = t_0$ and normal to this sphere.
	
	The initial data is required to satisfy the following bounds: for all $n \leq N_2$, we have
	\begin{equation*}
	\int_{\Sigma_{\tau_0}} \bigg(
	(1+r)^{-C_{[n]}\epsilon} |\overline{\slashed{\D}} \mathscr{Y}^n \phi|^2_{(A)} 
	+ (1+r)^{1-C_{[n]}\epsilon} |\slashed{\D}_L \mathscr{Y}^n \phi|^2_{(A)} 
	\bigg)\dVol_g \leq \epsilon^{2(N_2 + 3 - n)}
	\end{equation*}
	as well as the pointwise bounds
	\begin{equation*}
	\int_{\bar{S}_{t,r}} |\mathscr{Y}^n \phi|_{(A)} \dVol_{\mathbb{S}^2} \lesssim \epsilon^{2(N_2 + 3 - n)} (t + 1 - t_0)^{-1 + \frac{1}{2}C_{[n]}\epsilon}
	\end{equation*}
	Again, see \cite{Keir2018} for the definitions of the volume forms and the spheres $\bar{S}_{t,r}$.

	Then, if $N_2 \geq 8$ and $N_2 - 4 \geq N_1 \geq 4$, for all sufficiently small $\epsilon$ the system of wave equations $\tilde{\Box}_g \phi = F$ has a \emph{unique, global solution}, i.e.\ a unique solution in the region to the future of $\Sigma_{\tau_0}$. Furthermore, this solution will obey the pointwise bounds and $L^2$ bounds given in chapter 15 of \cite{Keir2018}, as well as the $L^2$ bounds\footnote{In fact, the system will obey all of these bounds with an additional factor of, say $1/2$ on the right hand side.} given in chapter 13.

\end{theorem}

\section{Proving theorem \ref{theorem main theorem}}
\label{section proving theorem}

The main theorem \ref{theorem main theorem} is identical to the theorem proved in our previous work \cite{Keir2018} (with a single level in the hierarchy), except for the fact that the stable boundedness property of the asymptotic system \ref{Condition boundedness and stability} replaces the hierarchical null condition. As such, rather than proving this theorem from scratch, we will simply show that each time the hierarchical null condition is used, condition \ref{Condition boundedness and stability} can be used as an alternative, leading to the same conclusions.

The hierarchical null condition is used only twice in the proof presented in our previous work \cite{Keir2018}. First, it is used when proving the energy estimates: if we are attempting to prove an estimate for a field at level $m$ in the hierarchy, we encounter error terms with ``bad weights'' in $r$, but which appear multiplying terms involving fields at lower levels in the hierarchy. For example, when attempting to prove the boundedness of the energy on a surface $\Sigma_{\tau_2}$ for a field $\phi_{(m)}$ at level $m$ in the hierarchy, we will encounter error terms of the form
\begin{equation*}
	\int_{\tau_1}^{\tau_2} \int_{\Sigma_\tau} \epsilon (1+r)^{-1 + C\epsilon} w_{(m)} |\partial \phi_{(m-1)}|^2 \dVol_{\Sigma_\tau} \upd \tau
\end{equation*}
where $w_{(m)} = (1+r)^{-C_{(m)}\delta}$ is the ``degenerate weight'' that we choose for the field $\phi_{(m)}$. In \cite{Keir2018} we dealt with these terms by choosing $w_{(m)} \gg w_{(m-1)}$, in other words, we chose the weights to be increasingly degenerate as we rise up the hierarchy.

Note that, at the bottom level of the hierarchy, we do not encounter this kind of error term. Instead, we encountered terms of the form
\begin{equation*}
	\int_{\tau_1}^{\tau_2} \int_{\Sigma_\tau} \epsilon (1+r)^{-1} w_{(0)} |\partial \phi|^2 \dVol_{\Sigma_\tau} \upd \tau
\end{equation*}
We deal with these terms by simply choosing $w_{(0)} = (1+r)^{-\delta}$: this generates bulk terms in the energy estimates that we can use to control these error terms. If we treat the fields in the problem that we are now considering as fields at the ``bottom level'', then the new semilinear terms that we encounter are exactly of this form. In other words, as long as we can recover the pointwise bounds
\begin{equation*}
	|\partial \phi_{(A)}| \lesssim \epsilon (1+r)^{-1}
\end{equation*}
for \emph{all} of the fields $\phi_{(A)}$, then the new error terms that we encounter in the energy estimate are of exactly the same form as error terms that we already encountered (and dealt with) in \cite{Keir2018}.

The other place in which the hierarchical null condition is used in \cite{Keir2018} is in recovering the sharp pointwise bounds. Since we use a degenerate energy, the usual approach of combining energy estimates with Sobolev embedding (on the spheres) will only yield the pointwise decay rate
\begin{equation*}
	|\partial \phi_{(m)}| \lesssim \epsilon (1+r)^{-1 + C_{(m)}\delta}
\end{equation*}
for fields at level $m$ in the hierarchy. However, in order to close our estimates, we need to upgrade these to the sharp decay rates
\begin{equation*}
\begin{split}
	|\partial \phi_{(m)}| &\lesssim \epsilon (1+r)^{-1 + C_{(m)}\epsilon}
	\\
	|\partial \phi_{(0)}| &\lesssim \epsilon (1+r)^{-1}
\end{split}
\end{equation*}
with corresponding estimates for fields after commuting\footnote{Note that, each time we commute, we can expect to lose another $C\epsilon$ of decay. Hence, if $\mathscr{Y}$ is a commutation operator, then the sharp decay rates are $|\partial \mathscr{Y}^n \phi_{(m)}| \lesssim \epsilon (1+r)^{-1 + C_{(n,m)}\epsilon}$.}.

In the present problem, all of the fields are treated as fields at the ``bottom level'' of the hierarchy. Thus we must recover the sharp decay rates
\begin{equation*}
	|\partial \phi_{(A)}| \lesssim \epsilon (1+r)^{-1}
\end{equation*}
for \emph{all} of the fields $\phi_{(A)}$. In \cite{Keir2018} these rates were recovered by using the hierarchical structure: the fields at the bottom level of the hierarchy obey wave equations with nonlinear terms that essentially satisfy the classical null condition. Using this fact, and integrating along integral curves of $L$, we recovered the sharp decay rates.

Since we are not assuming the hierarchical null structure, we have to find a different way to prove these sharp decay rates. This is where we can make use of the stability and boundedness condition \ref{Condition boundedness and stability}. We do this as follows: the wave equations \eqref{equation system of nonlinear reduced wave equations} can be written as
\begin{equation}
\label{equation asymptotic system with error 1}
	L(r\Lbar \phi_{(A)})
	=
	r^{-1} F^{(BC)}_{(A)} (r\Lbar \phi_{(B)}) (r\Lbar \phi_{(C)})
	+ \text{Err}
\end{equation}
where the error terms can be shown, using only the \emph{weak} decay estimates coming from the energy estimates and Sobolev embedding, to obey the bounds
\begin{equation*}
	|\text{Err}|
	\lesssim
	\mathcal{E}_0 (1+\tau)^{-\beta}(1+r)^{-\delta}
\end{equation*}
where the constant $\mathcal{E}_0$ depends only on a certain (weighted, higher order) energy of the initial data (see \cite{Keir2018} for the details).

Now, if we choose $s = \log(2+r)$ and we note that $L(r) = 1$, we find that equation  \eqref{equation asymptotic system with error 1} can be written as the following system of ODEs along the integral curves of $L$:
\begin{equation}
\label{equation asymptotic system with error 2}
	\frac{\partial}{\partial s}(r\Lbar \phi_{(A)})
	=
	F^{(BC)}_{(A)} (r\Lbar \phi_{(B)}) (r\Lbar \phi_{(C)})
	+ \text{Err}
\end{equation}
where the error term satisfies the bound
\begin{equation*}
	|\text{Err}|
	\lesssim
	\mathcal{E}_0 (1+\tau)^{-\beta} e^{-\delta s}
\end{equation*}
Hence, using condition \ref{Condition boundedness and stability}, we find that, if the size of the initial data $\mathcal{E}_0$ is sufficiently small, and if $r\Lbar\phi_{(A)}$ satisfies the bound
\begin{equation*}
	|r\Lbar\phi_{(A)}| \Big|_{r = r_0} \leq \mathcal{E}_0 (1+\tau)^{-\beta}
\end{equation*}
then the fields $r\Lbar\phi_{(A)}$ will satisfy the bound
\begin{equation*}
	|r\Lbar\phi_{(A)}| \lesssim \mathcal{E}_0 (1+\tau)^{-\beta}
\end{equation*}
in the region $r \geq r_0$. Combining this with the weak decay estimates in the region $r \leq r_0$, we can obtain the bounds
\begin{equation*}
	|\Lbar\phi_{(A)}| \lesssim \mathcal{E}_0 (1+\tau)^{-\beta}(1+r)^{-1}
\end{equation*}
In other words, given sufficiently small initial data, we can upgrade the weak decay estimates to obtain the sharp decay rates for all of the fields $\phi_{(A)}$.

Now, we also have to obtain sharp decay estimates for commuted fields. However, in this case, we can follow the computations in \cite{Keir2018} exactly. Indeed, note that, after commuting with an operator $\mathscr{Y}$, the new semilinear terms are of the form $(\partial \phi)(\partial \mathscr{Y} \phi)$ -- that is, they are \emph{linear} in the first derivatives of $\mathscr{Y} \phi$. Such terms appear in \cite{Keir2018}, and are easily handled\footnote{Note, however, that we will lose some decay: after commuting $n$ times, we will be able to recover the pointwise bound $|\partial \mathscr{Y}^n \phi| \lesssim \epsilon (1+\tau)^{-\beta}(1+r)^{-1+C_{(n)}\epsilon}$. However, we only need the sharp decay (that is, without the $\epsilon$ loss) for first derivatives in order to close the estimates.}.

\section{Asymptotic systems obeying condition \ref{Condition boundedness and stability}}
\label{section Euler equations}

We will now present a family of asymptotic systems obeying condition \ref{Condition boundedness and stability}. Members of this family satisfy condition \ref{Condition boundedness and stability} due to the presence of an underlying Hamiltonian structure, which in turn provides the necessary boundedness and stability criteria.

The most obvious way to make use of a Hamiltonian structure would be to examine asymptotic systems which correspond directly to Hamilton's equations for some dynamical system. However, this approach fails for the following reason: the Hamiltonian function should be at a local minimum when the values of all of the fields vanish, since this would guarantee the required stability and boundedness properties. On the other hand, the asymptotic system is a system of ODEs with quadratic nonlinear terms (and no linear, zero-th order terms), meaning that, if this system corresponds to Hamilton's equations, then the corresponding Hamiltonian should vanish \emph{cubically} at the origin.

Instead, we will focus on left-invariant Hamiltonians on Lie groups. In this case, from Hamilton's equations a related set of equations can be derived, called the \emph{Euler equations}. These equations will be quadratically nonlinear even when the Hamiltonian function itself vanishes quadratically at the identity. The following discussion closely follows \cite{Mienko1978}. Much of this discussion is technically unnecessary for our desired result -- we could, instead, have simply stated the family of asymptotic systems and shown that, due to the presence of a suitable conserved quantity, they obey condition \ref{Condition boundedness and stability}. However, it is interesting to understand how these systems arise from an underlying Hamiltonian structure.

Let $G$ be some finite-dimensional Lie group, and let $H$ be a non-negative function on the cotangent bundle of $G$:
\begin{equation*}
\begin{split}
	H : T^*(G) &\rightarrow \mathbb{R} \\
	H &\geq 0
\end{split}
\end{equation*}

There is a natural symplectic form $\omega$ on the cotangent bundle $T^*(G)$ given, in terms of local coordinates\footnote{Here the $x^a$ are coordinates for $G$ and the $p_a$ are coordinates for a point in the cotangent space at the point with coordinates $x^a$.} $(x^a,p_a)$ by $\omega = \upd x^a \wedge \upd p_a$. Using this, we can define a vector field $V_H$ on the cotangent bundle $T^*(G)$ by the requirement that, for all vector fields $Y$ on the cotangent bundle, we have
\begin{equation*}
	\omega(V_H, Y) = \upd H(Y) 
\end{equation*}
This vector field generates the \emph{Hamiltonian flow} associated with $H$. Hamilton's equations correspond to the flow generated by this vector field.

We define ``left translation'' as the following action of the Lie group $G$ on itself: given $g \in G$,
\begin{equation*}
\begin{split}
	L_g : G &\rightarrow G \\
	h &\mapsto gh
\end{split}
\end{equation*}
i.e.\ $L_g$ simply corresponds to multiplication on the left by $g$. We can lift this to an action on the cotangent bundle:
\begin{equation*}
\begin{split}
	L_g^* : T^*(G) &\rightarrow T^*(G) \\
	x &\mapsto L_g^*(x)
\end{split}
\end{equation*}

Now, suppose that the Hamiltonian function $H$ is also \emph{left-invariant}, that is, $H$ is invariant under the diffeomorphisms of the cotangent bundle generated by left translations: for all $x \in T^*(G)$ and for all $g \in G$ we have
\begin{equation*}
	H(x)
	=
	H\left(L_g^*(x) \right)
\end{equation*}
It is not hard to see that, in such a case, the Hamiltonian vector field $V_H$ is also left-invariant, that is, for all $g \in G$, and for all $x \in T^*(G)$, we have
\begin{equation*}
	dL^*_g (V_H)\big|_x = V_H \big|_{L^*_g(x)}
\end{equation*}
If we define $\pi$ as the canonical projection
\begin{equation*}
	\pi : T^*(G) \rightarrow G
\end{equation*}
then we can also define the projection from the cotangent bundle to the cotangent space at the origin, using left translation:
\begin{equation*}
\begin{split}
	l: T^*(G) &\rightarrow T^*_e(G)
	\\
	x &\mapsto L^*_{(\pi(x))^{-1}}(x)
\end{split}
\end{equation*}
We can also consider the map
\begin{equation*}
	\upd l: TT^*(G) \rightarrow T^*_e(G)
\end{equation*}
where we are implicitly using the canonical isomorphism $T_pT_e^*(G) \cong T_e^*(G)$ for $p \in T^*(G)$. Applying this to the vector $V_H$, we obtain a covector field $E_H$ on the cotangent space at the origin. This covector field depends, \emph{a priori}, on the point $x$ in the cotangent bundle $\in T^*(G)$ at which the vector field $V_H$ is evaluated, however, since $V_H$ is left-invariant, it in fact depends only on the value of $l(x) \in T^*_e(G)$. In other words, we can use $\upd l$ to translate the flow generated by $V_H$ on the cotangent bundle into a flow on the cotangent space at the identity, generated by the covector field $E_H$.

If we now introduce the bracket
\begin{equation*}
\begin{split}
	\{ \cdot, \cdot \} : T^*_e(G) \times T_e(G) &\rightarrow T^*(T_e(G))
	\\
	\{ \xi, X \}(Y) &= \xi \left( [X, Y] \right)
\end{split}
\end{equation*}
where $[\cdot, \cdot]$ is the Lie bracket of the group $G$, and we write $\bar{H}$ for the restriction of $H$ to $T^*_e(G)$ then it is shown in \cite{Mienko1978} that, setting $y = l(x) \in T_e^*(G)$, $E_H$ can be expressed as
\begin{equation*}
	E_H\big|_{y}
	=
	\{ y, \upd \bar{H}\big|_{y} \}
\end{equation*}
Note that, since $\bar{H}$ is a smooth function on $T_e^*(G)$, $\upd \bar{H}$ naturally takes values in the dual space to $T_e^*(G)$, that is, $T_e(G)$, via the following map: for all $X \in T^* T^*_e(G)$, we identify $X$ with $\tilde{X} \in T_e(G)$, where, for all $\xi \in T^*_e(G)$ we have
\begin{equation*}
	X (\xi) = -\xi (\tilde{X})
\end{equation*}
Hence, when we write $\upd \bar{H}$ in the expression for $E_H$ written above, we really mean the vector field in $T_e(G)$ corresponding to $\upd \bar{H}$ via the above correspondence.

The \emph{Euler equations} corresponding to the left-invariant Hamiltonian $H$ on the Lie group $G$ are given by the flow generated by $E_H$ on the cotangent space at the identity. In other words, the Euler equations are the following autonomous system of ODEs for the covector $y(s) \in T^*_e(G)$:
\begin{equation*}
	\frac{\upd}{\upd s} y(s) = E_H (y(s))
\end{equation*}

Now, note that $\bar{H}(y(s))$ is constant along the flow generated by the covector field $E_H$: we have
\begin{equation*}
\begin{split}
	\frac{\upd}{\upd s} \bar{H}(y(s))
	&=
	\frac{\upd}{\upd s} y(s) \cdot \upd \bar{H} \big|_{y(s)}
	\\
	&=
	\{ y(s), \upd \bar{H}\big|_{y(s)} \} (\upd \bar{H} \big|_{y(s)})
	\\
	&=y(s) \cdot \left[ \bar{H}\big|_{y(s)}\ , \ \bar{H}\big|_{y(s)} \right]
	\\
	&= 0
\end{split}
\end{equation*}	
by the antisymmetry of the Lie bracket.

Some of these expressions may be easier to understand if we introduce a basis $e^a$ for the cotangent space $T^*_e(G)$, and the dual basis $e_a$ for the tangent space $T_e(G)$ (so $e^a (e_b) = \delta_b^a$). Then, the components of $y$ satisfy the system of ODEs
\begin{equation*}
	\frac{\upd}{\upd s} y_a
	=
	C^c_{ba} y_c \frac{\partial \bar{H}}{\partial y_b} \bigg|_y
\end{equation*}
where $C^c_{ba}$ are the structure coefficients of the Lie algebra in this basis, i.e.
\begin{equation*}
	[e_b, e_a] = C^c_{ba} e_c
\end{equation*}
Then it is also clear that $\bar{H}$ is conserved by this flow: we have
\begin{equation*}
\begin{split}
	\frac{\upd}{\upd s} \bar{H}(y(s))
	&=
	\frac{\partial \bar{H}}{\partial y_a} \bigg|_y \frac{\upd y}{\upd s}
	\\
	&=
	C^c_{ba} y_c \left( \frac{\partial \bar{H}}{\partial y_b} \bigg|_y \right)  \left( \frac{\partial \bar{H}}{\partial y_a} \bigg|_y \right)
	\\
	&= 0
\end{split}
\end{equation*}
by the antisymmetry of the structure coefficients.

Now, suppose that $\bar{H}(y)$ is quadratic in $y$ and vanishes at the origin: that is, $\bar{H}(y) = \tilde{H}(y, y)$, where $\tilde{H}$ is linear in both of its arguments. Furthermore, suppose that the quadratic form $\tilde{H}$ is positive and non-degenerate, so that $\tilde{H}(y,y) \geq 0$ with equality if and only if $y = 0$. Then the equations for $y_a$ can be written as
\begin{equation*}
	\frac{\upd}{\upd s} y_a
	=
	C^c_{ba} \tilde{H}^{bd} y_c y_d 
\end{equation*}
where $\tilde{H}(x, y) = \tilde{H}^{ab} x_a y_b$.

Systems of this sort are systems of nonlinear ODEs with quadratic nonlinearities -- they can therefore arise as asymptotic systems associated with nonlinear systems of wave equations. We will now show that these asymptotic systems also obey condition \ref{Condition boundedness and stability}.

Consider the system
\begin{equation*}
	\frac{\upd}{\upd s} y_a
	=
	C^c_{ba} \tilde{H}^{bd} y_c y_d
	+ F_a(s)
\end{equation*}
where
\begin{equation*}
	|F_{(a)}(s)| \leq C\epsilon e^{-\delta s}
\end{equation*}
Let us consider the change in the value of the Hamiltonian $\bar{H}(y) = \tilde{H}(y,y)$ as $s$ evolves from $s = 0$ to $s = \infty$. An easy calculation shows that we have
\begin{equation*}
	\frac{\upd}{\upd s} \bar{H}(y(s))
	=
	2\tilde{H}^{ab} y_a F_b
\end{equation*}
Now, we can use $\tilde{H}$ to define an inner product on $T_e^*(G)$, via
\begin{equation*}
	\langle x , y \rangle := \tilde{H}(x,y)
\end{equation*}
then, writing $|y|^2 := \langle y, y \rangle$, we have\footnote{Note that, since the $y_a$'s lie in a finite dimensional vector space, all norms are equivalent.}
\begin{equation*}
\begin{split}
	\frac{\upd}{\upd s} |y|^2 = 2\langle y, F \rangle
	\\
	&\leq 2|y| |F|
\end{split}
\end{equation*}
and consequently, if $|y| \neq 0$,
\begin{equation*}
	\frac{\upd}{\upd s} |y| \leq |F| \lesssim \epsilon e^{-\delta s}
\end{equation*}
Integrating this to $s = \infty$ we find that
\begin{equation*}
	\lim_{s\rightarrow\infty} \Big| |y(s)| - |y(0)| \Big| \lesssim \epsilon \delta^{-1}
\end{equation*}
and so these systems obey condition \ref{Condition boundedness and stability}.

We obtain the following corollary of theorem \ref{theorem main theorem}:
\begin{corollary}
	Consider a system of nonlinear wave equations with an asymptotic system of the form
	\begin{equation*}
		\frac{\upd}{\upd s} \Phi_{(A)} = C^{(B)}_{(D)(A)} \tilde{H}^{(D)(E)} \Phi_{(B)}\Phi_{(E)}
	\end{equation*}
	where $C^{(B)}_{(D)(A)} = -C^{(B)}_{(A)(D)}$ are constants, as are $\tilde{H}^{(D)(E)}$, and where $\tilde{H}^{(D)(E)} y_{(D)} y_{(E)} \geq 0$ with equality if and only if $y = 0$. 
	
	Then this system of nonlinear wave equations admits global solutions for all sufficiently small initial data, where ``sufficiently small'' is defined in theorem \ref{theorem main theorem}.
\end{corollary}

\sloppy
\printbibliography

\end{document}